\documentclass[12pt]{article}
\usepackage{mathrsfs}
\usepackage{amsmath,amssymb}
\newcommand{\be}{\begin{equation}}
\newcommand{\ee}{\end{equation}}
\newcommand{\bal}{\begin{aligned}}
\newcommand{\eal}{\end{aligned}}
\newcommand{\bee}{\begin{equation*}}
\newcommand{\eee}{\end{equation*}}
\usepackage{color}

\def\cwedge{\bigcirc\kern-1.07em\wedge\ }

\newtheorem{thm}{Theorem}[section]

\newtheorem{prop}[thm]{Proposition}

\newtheorem{rem}{Remark}
\newtheorem{ex}{Example}
\newtheorem{defn}{Definition}
\numberwithin{equation}{section}

\begin{document}
\begin{center}

{\large {\bf  A generalization of a 4-dimensional Einstein
manifold}}
\end{center}
\footnotetext{\small{\it E-mail addresses}: {\bf
prettyfish@skku.edu} (Y. Euh), {\bf parkj@skku.edu} (J. H. Park),
{\bf sekigawa@math.sc.niigata-u.ac.jp} (K. Sekigawa).}

\begin{center}
Yunhee Euh${}^{\dag}$, JeongHyeong Park ${}^{\ddag}$ and Kouei
Sekigawa$^{\dag}$

\end{center}
{\small
\begin{center}
$~~~{}^{\dag}$Department of Mathematics,
    Niigata University,
    Niigata 950-2181, JAPAN\\
$~~~{}^{\ddag}$ Department of Mathematics,
    Sungkyunkwan University,
    Suwon 440-746, KOREA
\end{center}
}
\begin{abstract}

A weakly Einstein manifold is a generalization of a 4-dimensional
Einstein manifold, which is defined as an application of a curvature
identity derived from the generalized Gauss-Bonnet formula for a
4-dimensional compact oriented Riemannian manifold. In this paper,
we shall give a characterization of a weakly Einstein manifold.
\end{abstract}
\noindent {\it Mathematics Subsect Classification (2010)} : 53B20, 53C20 

\section{Introduction}

In the previous paper \cite{EPS}, we derived a curvature identity on
a 4-dimensional compact oriented Riemannian manifold from the
generalized Gauss-Bonnet formula, and further gave a direct proof of
the fact that the curvature identity holds on any 4-dimensional
Riemannian manifold which is not necessarily compact. Consequently,
we proved that the following curvature identity holds on any
$4$-dimensional Riemannian manifold $M=(M,g)$:
   \be\label{main equation}
    \check{R}-2\check{\rho}-L\rho+\tau\rho-\frac{1}{4}(|R|^2-4|\rho|^2+\tau^2)g=0.
    \ee
    Here,
     \bee
     \begin{gathered}
     \check{R}:\check{R}_{ij}=R_{abci}R^{abc}_{~~~j},\qquad
     \check{\rho}:\check{\rho}_{ij}=\rho_{ai}\rho^{a}_{~j},\\
     L:(L\rho)_{ij}=2R_{iabj}\rho^{ab},
     \end{gathered}
    \eee
     where $R$, $\rho$ and $\tau$ are the curvature tensor, the Ricci tensor and   the scalar curvature of $M$, respectively.

From \eqref{main equation}, we may easily check that any
4-dimensional Einstein manifold $M=(M,g)$ satisfies the condition
    \be\label{eq:super-Einstein}
    {R^{abc}}_{i}R_{abcj}=\frac{1}{4}|R|^2g_{ij}.
    \ee

In \cite{EPS}, we defined a weakly Einstein manifold based on the
above, i.e. Riemannian manifold $M=(M,g)$ satisfying the condition
\eqref{eq:super-Einstein} (with $|R|^2$ not necessarily constant).
By the definition, we see immediately that a weakly Einstein
manifold in dimension 4 is a generalization of a 4-dimensional
Einstein manifold (see Examples 4 and 5 in \textsection 3). We may
also remark that a weakly Einstein manifold is not necessarily
Einstein. As a characterization of a 4-dimensional Einstein
manifold, the following theorem is well-known.\vskip0.2cm

\noindent{\bf Theorem A} (\cite{S-T}) {\it A 4-dimensional
Riemannian manifold $M=(M,g)$ is Einstein if and only if there
exists a Singer-Thorpe basis of $T_pM$ at each point $p\in
M$.}\vskip0.2cm

The main purpose of the present paper is to give a generalization of
Theorem A. Namely, we shall prove the following:\vskip0.2cm

\noindent {\bf Theorem B} {\it A 4-dimensional Riemannian manifold
$M=(M,g)$ is weakly Einstein if and only if there exists a
generalized Singer-Thorpe basis of $T_pM$ at each point $p\in
M$.}\vskip0.2cm

In \textsection 2, we shall prepare some fundamental terminologies
and notational conventions for the forthcoming arguments. In
\textsection 3, we shall give a proof of Theorem B.

\section{Preliminaries}
Let $M=(M,g)$ be a 4-dimensional Riemannian manifold and
$\mathfrak{X}(M)$ be the Lie algebra of all smooth vector fields on
$M$. We denote the Levi-Civita connection, the curvature tensor, the
Ricci tensor and the scalar curvature of $M$ by $\nabla$, $R$,
$\rho$ and $\tau$, respectively. We assume that the curvature tensor
$R$ is defined by $R(X,Y)Z=[\nabla_X,\nabla_Y]Z-\nabla_{[X,Y]}Z$ for
$X$, $Y$, $Z\in\mathfrak{X}(M)$. Further, we denote the Ricci
transformation by $Q$ given by $\rho(X,Y)=g(QX,Y)$ for $X$, $Y\in
\mathfrak{X}(M)$. Then, we may easily check that $Q$ is symmetric
with respect to the metric $g$, namely, $g(QX,Y)=g(X,QY)$ for $X$,
$Y\in \mathfrak{X}(M)$. Now, we may rewrite the curvature identity
\eqref{main equation} as follows:
    \be\label{eq:id}
    \bal
    &\sum_{a,b,c}R_{abci}R_{abcj}-2\sum_{a}\rho_{ai}\rho_{aj}-2\sum_{a,b}\rho_{ab}R_{iabj}\\
    &\qquad+\tau\rho_{ij}-\frac{1}{4}\big(|R|^2-4|\rho|^2+\tau^2\big)\delta_{ij}=0,
    \eal
    \ee
 with respect to an orthonormal basis $\{e_{i}\}$ ($1 \leq i
          \leq 4$) of $T_{p}M$ at any point $p \in M$, where
          $R_{ijkl}=g(R(e_i,e_j)e_k, e_l)$,
          $\rho_{ij}=\rho(e_i,e_j)$.

We here introduce some special kinds of orthonormal basis of $T_pM$
($p\in M$) and explain their intermediate relationships. We assume
that an orthonormal basis $\{e_i\}$ ($1\leq i\leq 4$) of $T_pM$ is
simultaneously a Ricci eigenbasis and Chern basis \cite{Chern, Kl,
K-P} satisfying
   \be\label{eq:Chern_basis}
    R_{1213}=R_{1214}=R_{1223}=R_{1224}=R_{1314}=R_{1323}=0.
    \ee
Then, we have further
  \be\label{eq:Chern_basis3}
    R_{2434}=R_{2334}=R_{1434}=R_{1334}=R_{2324}=R_{1424}=0.
    \ee
Thus, from \eqref{eq:Chern_basis} and \eqref{eq:Chern_basis3}, we
have
    \be\label{eq:general_S-T}
    R_{ijjk}=0~~(i\ne k),\quad (1\leq i,j,k\leq 4).
    \ee
Conversely, if \eqref{eq:general_S-T} holds with respect to an
orthonormal basis $\{e_i\}$ of $T_pM$, then we see that the basis
$\{e_i\}$ is a Ricci eigenbasis and a Chern basis at the same time.

The following example shows that a Ricci eigenbasis is not
necessarily always a Chern basis.
\begin{ex}\label{ex:5}
{\rm Let $\mathfrak{g}=\text{span}_{\mathbb{R}}\{e_1,e_2,e_3,e_4\}$
be a 4-dimensional real Lie algebra equipped with the following Lie
bracket operation:
    \be\label{eq:LieBracket_ex5}
    \bal
    &[e_1,e_2]=2e_2,\qquad[e_1,e_3]=-e_3,\qquad[e_1,e_4]=2e_3-e_4,\\
    &[e_2,e_3]=0,\quad\qquad[e_2,e_4]=0,\quad\qquad[e_3,e_4]=0,
    \eal
    \ee
and $<,>$ the inner product on $\mathfrak{g}$ given by
$<e_i,e_j>=\delta_{ij}$. Let $G$ be a connected and simply connected
solvable Lie group with the Lie algebra $\mathfrak{g}$ of $G$ and
$g$ the $G$-invariant Riemannian metric on $G$ determined by $<,>$.
We set $\nabla_{e_i}e_j=\sum_{k=1}^4\Gamma_{ijk}e_k$, $(1\leq
i,j\leq 4)$. Then, we get
    \be\label{def:Gamma}
    \Gamma_{ijk}=-\Gamma_{ikj}
    \ee
and further, from \eqref{eq:LieBracket_ex5}, we obtain
    \be\label{eq:Gamma_ex5}
    \bal
    &\Gamma_{134}=-1,\quad \Gamma_{212}=-2,\quad \Gamma_{313}=1,\\
    &\Gamma_{314}=-1,\quad\Gamma_{413}=-1,\quad\Gamma_{414}=1,
    \eal
    \ee
and otherwise being zero up to sign. From \eqref{def:Gamma}
 and \eqref{eq:Gamma_ex5}, by direct calculations, we
have
    \be
    \bal
    &R_{1212}=4,~~\quad R_{1414}=4,\\
    &R_{2323}=-2,\quad R_{2424}=-2,\\
    &R_{1314}=-2,\quad R_{2324}=2,
    \eal
    \ee
and otherwise being zero up to sign. Then, we have the Ricci
eigenvalues as $\lambda_1=-8$, $\lambda_2=0$, $\lambda_3=2$,
$\lambda_4=-2$.}
\end{ex}
Now, we recall the definition of a Singer-Thorpe basis. An
orthonormal basis $\{e_i\}$ of $T_pM$ ($p\in M$) is called a {\it
Singer-Thorpe basis} at $T_pM$ if the basis $\{e_i\}$ satisfies
\eqref{eq:general_S-T} and
    \be\label{eq:Siner-Thorpe}
    R_{1212}=R_{3434},\quad R_{1313}=R_{2424},\quad
    R_{1414}=R_{2323}.
    \ee
We here give a generalization of the Singer-Thorpe basis.
\begin{defn}
Let $M=(M,g)$ be a 4-dimensional Riemannian manifold and $\{e_i\}$
be an orthonormal basis of $T_pM$ at $p\in M$. If the basis
$\{e_i\}$ satisfies \eqref{eq:general_S-T} and
    \be\label{eq:general-Singer}
    {R_{1212}}^2={R_{3434}}^2,\quad{R_{1313}}^2={R_{2424}}^2,\quad{R_{1414}}^2={R_{2323}}^2,
    \ee
then the orthonormal basis is called a {\it generalized
Singer-Thorpe basis} of $T_pM$.
\end{defn}

\section{Proof of Theorem B}
First, we shall prove the following proposition which gives a
necessary condition for a 4-dimensional Riemannian manifold to be
weakly Einstein.
\begin{prop}\label{th:4.1}
Let $M=(M,g)$ be a weakly Einstein manifold and $\{e_i\}$ $(1\leq
i\leq4)$
  an orthonormal Ricci eigenbasis of $T_pM$
corresponding to the eigenvalues $\lambda_i$ $(1\leq i\leq4)$ at any
point $p\in M$. Then, we see that the curvature condition
    \be\label{eq:from(5.1)}
    {R_{1212}}^2={R_{3434}}^2,\quad
    {R_{1313}}^2={R_{2424}}^2,\quad
    {R_{1414}}^2={R_{2323}}^2
    \ee
holds and also the following cases $(1)\sim(4)$ never
occur:\\
$(1)\quad \lambda_1=\lambda_2=\lambda_3(\ne0),~\lambda_4=0,$\\
$(2)\quad \lambda_1=\lambda_2=\lambda_4(\ne0),~\lambda_3=0,$\\
$(3)\quad \lambda_1=\lambda_3=\lambda_4(\ne0),~\lambda_2=0,$\\
$(4)\quad \lambda_2=\lambda_3=\lambda_4(\ne0),~\lambda_1=0.$\\
Especially, if $M$ is Einstein, then
           \bee
           R_{1212} = R_{3434},\quad　R_{1313} = R_{2424}, \quad R_{1414}= R_{2323}
           \eee
　　　holds for any orthonormal basis $\{e_{i} \}$ of $T_{p}M$.
\end{prop}
{\it Proof.} Let $M=(M,g)$ be a weakly Einstein manifold and $p$ any
point of $M$ and $\{e_i$\} ($1\leq i\leq4$) an orthonormal Ricci
eigenbasis of $T_pM$ corresponding to the Ricci eigenvalues
$\lambda_i$ $(1\leq i\leq4)$ at $p$, namely, satisfying the
following condition
    \be\label{def:basis_Ric}
    Qe_i=\lambda_ie_i\quad(1\leq i\leq4).
    \ee
Then, from \eqref{def:basis_Ric}, we get
    \be\label{eq:norm_curv_ric_eigen}
    \bal
    |R|^2
    =&4\{{R_{1212}}^2+{R_{1313}}^2+{R_{1414}}^2+{R_{2323}}^2+{R_{2424}}^2+{R_{3434}}^2\\
    &+4{R_{1213}}^2+4{R_{1214}}^2+4{R_{1223}}^2+4{R_{1224}}^2+4{R_{1314}}^2+4{R_{1323}}^2\\
    &+2{R_{1234}}^2+2{R_{1342}}^2+2{R_{1423}}^2\}.
    \eal
    \ee
On the other hand, setting $i=j=1$ in the left hand side of
\eqref{eq:super-Einstein}, we get
    \be\label{eq:i=1,j=2in_as_lhs}
    \bal
    \sum_{a,b,c}{R_{abc1}}^2=&2\{{R_{1212}}^2+{R_{1313}}^2+{R_{1414}}^2+{R_{1234}}^2+{R_{1342}}^2+{R_{1423}}^2\\
    +&2({R_{1213}}^2+{R_{1214}}^2+{R_{1223}}^2+{R_{1224}}^2+{R_{1314}}^2+{R_{1323}}^2)\}.
    \eal
    \ee
From \eqref{eq:norm_curv_ric_eigen}, \eqref{eq:i=1,j=2in_as_lhs},
and taking account of \eqref{eq:super-Einstein}, we have the
following equality
     \be\label{eq:i=j=1in_as}
    \bal
   {R_{1212}}^2+{R_{1313}}^2+{R_{1414}}^2-{R_{2323}}^2-{R_{2424}}^2-{R_{3434}}^2=0.
    \eal
    \ee
Similarly, we get
    \be\label{eq:i=j=2in_as}
    \bal
    {R_{1212}}^2+{R_{2323}}^2+{R_{2424}}^2-{R_{1313}}^2-{R_{1414}}^2-{R_{3434}}^2=0,
    \eal
    \ee
    \be\label{eq:i=j=3in_as}
    \bal
     {R_{1212}}^2+{R_{1414}}^2+{R_{2424}}^2-{R_{1313}}^2-{R_{2323}}^2-{R_{3434}}^2=0,
    \eal
    \ee
    \be\label{eq:i=j=4in_as}
    \bal
    {R_{1212}}^2+{R_{1313}}^2+{R_{2323}}^2-{R_{1414}}^2-{R_{2424}}^2-{R_{3434}}^2=0.
    \eal
    \ee
From \eqref{eq:i=j=1in_as} and  \eqref{eq:i=j=2in_as}, we have
    \be\label{eq:R_1212}
    {R_{1212}}^2-{R_{3434}}^2=0.
    \ee
Similarly, from \eqref{eq:i=j=1in_as} and \eqref{eq:i=j=3in_as}, we
have
    \be\label{eq:R_1313}
    {R_{1313}}^2-{R_{2424}}^2=0.
    \ee
From \eqref{eq:i=j=1in_as} and  \eqref{eq:i=j=4in_as}, we have
    \be\label{eq:R_1414}
    {R_{1414}}^2-{R_{2323}}^2=0.
    \ee
Thus, from \eqref{eq:R_1212}$\sim$\eqref{eq:R_1414}, we have
\eqref{eq:from(5.1)}.

 Next, from \eqref{eq:from(5.1)}, we see that
the following
eight cases can be taken into consideration;\\
Case (i) $R_{1212}=R_{3434},~R_{1313}=R_{2424},~R_{1414}=R_{2323}$.\\
Then, $\lambda_1-\lambda_2=0,\quad \lambda_1-\lambda_3=0,\quad
\lambda_1-\lambda_4=0$, and hence,
$\lambda_1=\lambda_2=\lambda_3=\lambda_4$.\\
Case (ii) $R_{1212}=-R_{3434},~R_{1313}=R_{2424},~R_{1414}=R_{2323}$.\\
Then, we get also
    \bee
    \bal
    &\lambda_1-\lambda_2=0,\qquad\qquad\qquad\lambda_1-\lambda_3=-2R_{1212},\\
    &\lambda_1-\lambda_4=-2R_{1212},\qquad\quad\lambda_2-\lambda_3=-2R_{1212},\\
    &\lambda_2-\lambda_4=-2R_{1212},\qquad\quad\lambda_3-\lambda_4=0,
    \eal
    \eee
and hence, $\lambda_1=\lambda_2,~ \lambda_3=\lambda_4$.\\
Case (iii) $R_{1212}=R_{3434},~R_{1313}=-R_{2424},~R_{1414}=R_{2323}$.\\
Then, we get
    \bee
    \bal
    &\lambda_1-\lambda_2=-2R_{1313},\qquad\lambda_1-\lambda_3=0,\\
    &\lambda_1-\lambda_4=-2R_{1313},\qquad\lambda_2-\lambda_3=2R_{1313},\\
    &\lambda_2-\lambda_4=0,\qquad\qquad\quad\text{ }\lambda_3-\lambda_4=-2R_{1313},
    \eal
    \eee
and hence, $\lambda_1=\lambda_3,~ \lambda_2=\lambda_4$.\\
Case (iv) $R_{1212}=R_{3434},~R_{1313}=R_{2424},~R_{1414}=-R_{2323}$.\\
Then, we get
    \bee
    \bal
    &\lambda_1-\lambda_2=-2R_{1414},\qquad\lambda_1-\lambda_3=-2R_{1414},\\
    &\lambda_1-\lambda_4=0,\qquad\qquad\quad\lambda_2-\lambda_3=0,\\
    &\lambda_2-\lambda_4=2R_{1414},\text{ }\quad\quad\text{ }\lambda_3-\lambda_4=2R_{1414},
    \eal
    \eee
and hence, $\lambda_1=\lambda_4,~ \lambda_2=\lambda_3$.\\
Case (v) $R_{1212}=R_{3434},~R_{1313}=-R_{2424},~R_{1414}=-R_{2323}$. \\
Then, we get
    \bee
    \bal
    &\lambda_1+\lambda_2=-2R_{1212},\qquad\lambda_1-\lambda_3=-2R_{1414},\\
    &\lambda_1-\lambda_4=-2R_{1313},\qquad\lambda_2-\lambda_3=2R_{1313},\\
    &\lambda_2-\lambda_4=2R_{1414},\qquad\quad\lambda_3+\lambda_4=-2R_{1212},
    \eal
    \eee
and hence, $\lambda_1+\lambda_2=\lambda_3+\lambda_4$.\\
Case (vi)
$R_{1212}=-R_{3434},~R_{1313}=R_{2424},~R_{1414}=-R_{2323}$.\\
Then, we get
    \bee
    \bal
    &\lambda_1-\lambda_2=-2R_{1414},\qquad\lambda_1+\lambda_3=-2R_{1313},\\
    &\lambda_1-\lambda_4=-2R_{1212},\qquad\lambda_2-\lambda_3=-2R_{1212},\\
    &\lambda_2+\lambda_4=-2R_{1313},\qquad\lambda_3-\lambda_4=2R_{1414},
    \eal
    \eee
and hence, $\lambda_1+\lambda_3=\lambda_2+\lambda_4$.\\
Case (vii)
$R_{1212}=-R_{3434},~R_{1313}=-R_{2424},~R_{1414}=R_{2323}$.\\
Then, we get
    \bee
    \bal
    &\lambda_1-\lambda_2=-2R_{1313},\qquad\lambda_1-\lambda_3=-2R_{1212},\\
    &\lambda_1+\lambda_4=-2R_{1414},\qquad\lambda_2+\lambda_3=-2R_{1414},\\
    &\lambda_2-\lambda_4=-2R_{1212},\qquad\lambda_3-\lambda_4=-2R_{1313},
    \eal
    \eee
and hence, $\lambda_1+\lambda_4=\lambda_2+\lambda_3$.\\
Case (viii)
$R_{1212}=-R_{3434},~R_{1313}=-R_{2424},~R_{1414}=-R_{2323}.$\\
Then, we get
    \bee
    \bal
    &\lambda_1+\lambda_2=-2R_{1212},\qquad\lambda_1+\lambda_3=-2R_{1313},\\
    &\lambda_1+\lambda_4=-2R_{1414},\qquad\lambda_2+\lambda_3=2R_{1414},\\
    &\lambda_2+\lambda_4=2R_{1313},\qquad\quad\lambda_3+\lambda_4=2R_{1212},
    \eal
    \eee
and hence, $\lambda_1+\lambda_2+\lambda_3+\lambda_4=0$ (i.e.,
$\tau=0$).\\
Thus, from the above arguments in Cases (i)$\sim$(viii), we see that
the cases (1)$\sim$(4) in Proposition \ref{th:4.1} do not occur.
\hfill$\square$\medskip
\begin{rem}
{\rm In the proof of Proposition \ref{th:4.1}, we may note that
Cases (ii) to (iv) (also for Cases (v) to (vii), respectively) are
all essentially equivalent.}
\end{rem}
 The following examples illustrate Proposition \ref{th:4.1}.
Then, from the examples we can easily check that $M$ is not a weakly
Einstein manifold.
\begin{ex}\label{ex:1}
{\rm Let $M$ be a Riemannian product manifold of 2-dimensional
Riemannian manifolds of constant Gaussian curvatures $c_1$ and $c_2$
satisfying $c_1^2\ne c_2^2$. Then this implies that $M$ is not a
weakly Einstein manifold.}
\end{ex}
\begin{ex}\label{ex:2}
{\rm Let $M=(M,g)$ be a Riemannian product manifold of a
3-dimensional space of constant sectional curvature $c(\ne 0)$ and a
real line $\mathbb{R}$. From Proposition \ref{th:4.1}, we see that
$M$ is not a weakly Einstein manifold.}
\end{ex}

\begin{rem}
{\rm Based on Proposition \ref{th:4.1} and the related Examples
\ref{ex:1} and \ref{ex:2}, it may be seen that the statement ``for
any 4-dimensional Riemannian manifold one always gets
\eqref{eq:super-Einstein}" (\cite{Be}, pp. 165), is incorrect.}
\end{rem}

The following examples show that a weakly Einstein manifold is not
necessarily Einstein.
\begin{ex}\label{ex:3}
{\rm (\cite{EPS}) Let $M$ be a Riemannian product manifold of
2-dimensional Riemannian manifolds $M_1(c)$ and $M_2(-c)$ of
constant Gaussian curvatures $c$ and $-c$ $(c\ne0)$, respectively.
Then we can easily check that $M$ is not Einstein.
 We can also easily check that M satisfies \eqref{eq:super-Einstein}, thus
 $M$ is weakly Einstein. Further, $M$ belongs to Cases (ii), (vi), (vii)
 and (viii).}
\end{ex}

\begin{ex}\label{ex:4}
{\rm Let $\mathfrak{g}=\text{span}_{\mathbb{R}}\{e_1,e_2,e_3,e_4\}$
be a 4-dimensional real Lie algebra equipped with the following Lie
bracket operation:
    \be\label{eq:LieBracket}
    \bal
    &[e_1,e_2]=ae_2,\qquad[e_1,e_3]=-ae_3-be_4,\qquad[e_1,e_4]=be_3-ae_4,\\
    &[e_2,e_3]=0,\quad\qquad[e_2,e_4]=0,\qquad\qquad\qquad[e_3,e_4]=0,
    \eal
    \ee
where $a(\ne0)$, $b$ are constant. We define an inner product $<,>$
on $\mathfrak{g}$ by \\$<e_i,e_j>=\delta_{ij}$.  Let $G$ be a
connected and simply connected solvable Lie group with the Lie
algebra $\mathfrak{g}$ of $G$ and $g$ the $G$-invariant Riemannian
metric on $G$ determined by $<,>$.  From \eqref{eq:LieBracket},
    \be\label{eq:Gamma}
    \Gamma_{134}=-b,\quad\Gamma_{212}=-a,\quad
    \Gamma_{313}=a,\quad\Gamma_{414}=a,
    \ee
and otherwise being zero up to sign. From \eqref{def:Gamma} and
\eqref{eq:Gamma}, by direct calculations, we have
    \be
    \begin{gathered}
    R_{1212}=a^2,~~\quad R_{1313}=a^2,~~\quad R_{1414}=a^2,\\
    R_{2323}=-a^2,\quad R_{2424}=-a^2,\quad R_{3434}=a^2,
    \end{gathered}
    \ee
and otherwise being zero up to sign. From this, we can easily check
that $M$ is not Einstein since the Ricci curvature components
satisfy  $\rho_{11}= -3a^2 $ but $\rho_{22}= a^2$. We also can
easily check that M satisfies \eqref{eq:super-Einstein}, thus $M$ is
weakly Einstein. Then, we see that $(G,g)$ belongs to Case (v).}
\end{ex}

\begin{rem}
{\rm Jensen \cite{Je} proved that a 4-dimensional homogeneous
Einstein manifold is locally symmetric. We may easily check that
Example \ref{ex:4} is homogeneous but not locally symmetric. Thus,
Example \ref{ex:4} shows that Jensen's result does not necessarily
hold for weakly Einstein manifolds in general.}
\end{rem}

In the remainder of this section, we shall give a proof of Theorem
B.

Necessity:  From Proposition \ref{th:4.1}, it suffices to prove that
there exists an orthonormal Ricci eigenbasis $\{e_i\}$ of $T_p M$ at
each point  $p\in M$ which satisfies \eqref{eq:general_S-T}.
 Let $M=(M,g)$ be a weakly Einstein manifold.
Then, from \eqref{eq:id}, we have also the following equality
    \be\label{eq:id+assumption}
    2\sum_a\rho_{ai}\rho_{aj}+2\sum_{a,b}\rho_{ab}R_{iabj}-\tau\rho_{ij}-|\rho|^2\delta_{ij}+\frac{\tau^2}{4}\delta_{ij}=0.
    \ee
Setting $i=j=1$ in \eqref{eq:id+assumption}, we get
    \be\label{eq:i=j=1}
    2\lambda_1^2+2\sum_{i}\lambda_iR_{1ii1}-\Big(\sum_{i}\lambda_i\Big)\lambda_1-\sum_{i}\lambda_i^2+\frac14\Big(\sum_{i}\lambda_i\Big)^2=0.
    \ee
Similarly, we get
    \be\label{eq:i=j=3}
    \bal
    &2\lambda_2^2+2\sum_{i}\lambda_iR_{2ii2}-\Big(\sum_{i}\lambda_i\Big)\lambda_2-\sum_{i}\lambda_i^2+\frac14\Big(\sum_{i}\lambda_i\Big)^2=0,\\
    &2\lambda_3^2+2\sum_{i}\lambda_iR_{3ii3}-\Big(\sum_{i}\lambda_i\Big)\lambda_3-\sum_{i}\lambda_i^2+\frac14\Big(\sum_{i}\lambda_i\Big)^2=0,\\
    &2\lambda_4^2+2\sum_{i}\lambda_iR_{4ii4}-\Big(\sum_{i}\lambda_i\Big)\lambda_4-\sum_{i}\lambda_i^2+\frac14\Big(\sum_{i}\lambda_i\Big)^2=0.
    \eal
    \ee
 Further, setting $i=1$, $j=2$ in \eqref{eq:id+assumption}, we get the
following
    \be\label{eq:i=1,j=2}
    (\lambda_3-\lambda_4)R_{1323}=0.
    \ee
Similarly, we get
    \be\label{eq:i=1,j=3}
    \bal
    &(\lambda_2-\lambda_4)R_{1223}=0,\qquad(\lambda_2-\lambda_3)R_{1224}=0,\qquad(\lambda_1-\lambda_4)R_{1213}=0,\\
    &(\lambda_1-\lambda_3)R_{1214}=0,\qquad(\lambda_1-\lambda_2)R_{1314}=0.
    \eal
    \ee
Then, the following cases are considerable:

{\bf Case I} $~~\lambda_1=\lambda_2=\lambda_3=\lambda_4$.

{\bf Case II-1}
$~\lambda_1=\lambda_2(\equiv\lambda),\quad\lambda_3\ne\lambda_4,\quad(\lambda_3,\lambda_4\ne\lambda)$.

{\bf Case II-2}
$~\lambda_1=\lambda_3(\equiv\lambda),\quad\lambda_2\ne\lambda_4,\quad(\lambda_2,\lambda_4\ne\lambda)$.

{\bf Case II-3}
$~\lambda_1=\lambda_4(\equiv\lambda),\quad\lambda_2\ne\lambda_3,\quad(\lambda_2,\lambda_3\ne\lambda)$.

{\bf Case II-4}
$~\lambda_2=\lambda_3(\equiv\lambda),\quad\lambda_1\ne\lambda_4,\quad(\lambda_1,\lambda_4\ne\lambda)$.

{\bf Case II-5}
$~\lambda_2=\lambda_4(\equiv\lambda),\quad\lambda_1\ne\lambda_3,\quad(\lambda_1,\lambda_3\ne\lambda)$.

{\bf Case II-6}
$~\lambda_3=\lambda_4(\equiv\lambda),\quad\lambda_1\ne\lambda_2,\quad(\lambda_1,\lambda_2\ne\lambda)$.

{\bf Case III-1}
$~\lambda_1=\lambda_2(\equiv\lambda),\quad\lambda_3=\lambda_4(\equiv\mu),\quad(\lambda\ne\mu)$.

{\bf Case III-2}
$~\lambda_1=\lambda_3(\equiv\lambda),\quad\lambda_2=\lambda_4(\equiv\mu),\quad(\lambda\ne\mu)$.

{\bf Case III-3}
$~\lambda_1=\lambda_4(\equiv\lambda),\quad\lambda_2=\lambda_3(\equiv\mu),\quad(\lambda\ne\mu)$.

{\bf Case IV-1}
$~\lambda_1=\lambda_2=\lambda_3(\equiv\lambda),\quad\lambda_4\ne\lambda.$

{\bf Case IV-2}
$~\lambda_1=\lambda_2=\lambda_4(\equiv\lambda),\quad\lambda_3\ne\lambda.$

{\bf Case IV-3}
$~\lambda_1=\lambda_3=\lambda_4(\equiv\lambda),\quad\lambda_2\ne\lambda.$

{\bf Case IV-4}
$~\lambda_2=\lambda_3=\lambda_4(\equiv\lambda),\quad\lambda_1\ne\lambda.$

{\bf Case V} $~~\lambda_i\ne\lambda_j,\quad(i\ne j).$

\noindent {\bf Case I.} The existence of a generalized Singer-Thorpe
basis follows immediately from the construction of a Singer-Thorpe
basis.

\noindent {\bf Case V.} Then, from \eqref{eq:i=1,j=2} and
\eqref{eq:i=1,j=3}, we may immediately choose a generalized
Singer-Thorpe basis.

\noindent{\bf Case II-1.} Then, it suffices to consider Cases (v)
and (viii). First, we deal with Case (v). From \eqref{eq:i=1,j=2}
and \eqref{eq:i=1,j=3}, taking account of the equalities in Case
(v), we have
    \be\label{eq:case_ii_1_other}
    \bal
    &R_{1323}=0,~R_{1223}=0,~R_{1224}=0,\\
    &R_{1213}=0,~R_{1214}=0,~
    R_{1313}=R_{2323}.
    \eal
    \ee
Here, we note that all of the relations in
\eqref{eq:case_ii_1_other} and Case (v) are preserved under the
changes of the orthonormal basis satisfying the conditions of
 Case II-1. We denote the 2-dimensional subspace of $T_pM$
spanned $\{e_1, e_2\}$ by $V$. For any non-zero vector $x\in V$, we
denote by $x^{\perp}$ the vector in $V$ such that $|x^{\perp}|=|x|$,
$g(x,x^{\perp})=0$, and the ordered pair $\{x,x^{\perp}\}$ and
$\{e_1,e_2\}$ determine the same orientation on $V$. We define a
unit vector $e\in V$ by
    \be\label{eq:case_ii_1}
    R(e,e_3,e^{\perp},e_4)=\max_{x\in
    V,~|x|=1}R(x,e_3,x^{\perp},e_4).
    \ee
We set $e'_1=e$, $e'_2=e^{\perp}$, $e'_3=e_3$, $e'_4=e_4$ and define
a function $\phi(t)$ by
    \be\label{eq:case_ii_1_phi}
    \phi(t)=R(\cos t e'_1+\sin t e'_2,e'_3,-\sin t e'_1+\cos t
    e'_2,e'_4).
    \ee
Then, from \eqref{eq:case_ii_1} and \eqref{eq:case_ii_1_phi}, we
have $\phi'(0)=0,$ and hence,
    \be\label{eq:case_ii_1_1314}
    0=-R'_{1314}+R'_{2324}=-2R'_{1314}~(\text{and hence, }
    R'_{2324}=0),
    \ee
where $R'_{ijkl}=R(e'_i,e'_j,e'_k,e'_l)$, $1\leq i,j,k,l \leq 4.$
Then together with \eqref{eq:case_ii_1_1314}, the respective
equalities in \eqref{eq:case_ii_1_other} and Case (v) corresponding
to the orthonormal basis $\{e'_{i}\}$, we see that the orthonormal
basis $\{e'_i\}$ is a generalized Singer-Thorpe basis. Similarly, we
may also choose a generalized Singer-Thorpe basis for Case (viii).
Further, we may also choose a generalized Singer-Thorpe basis for
Cases II-2$\sim$II-6.

\noindent{\bf Case III-1.} Then it suffices to consider Cases (ii),
(vi), (vii), (viii). First, we consider Case (ii). Then, from
\eqref{eq:i=1,j=2} and \eqref{eq:i=1,j=3}, we have
    \be\label{eq:case_iii_1_other}
    \begin{gathered}
    R_{1223}=0,~R_{1224}=0,~R_{1213}=0,~R_{1214}=0.\\
    \end{gathered}
    \ee
Here, we may note that each of the relations in
\eqref{eq:case_iii_1_other} and Case (ii) is preserved under the
changes of the orthonormal basis satisfying the conditions of
 Case III-1. Let $V$ be a
2-dimensional subspace of $T_pM$ spanned by $\{e_1,e_2\}$ and
$V^{\perp}$ be the orthogonal complement of $V$ in $T_pM$. Then
$V^{\perp}$ is spanned by $\{e_3,e_4\}$. We define $e'_1\in V$ and
$e'_3\in V^{\perp}$ by
    \be\label{eq:case_iii_1}
    R(e'_1,e'_3,e'_1,e'_3)=\underset{|x|=|y|=1}{\max_{x\in V,~y\in
    V^{\perp}}}R(x,y,x,y).
    \ee
Further, we choose unit vectors $e'_2\in V$ and $e'_4\in V^{\perp}$
in such a way that $\{e_1,e_2\}$ and $\{e'_1,e'_2\}$ ($\{e_3,e_4\}$
and $\{e'_3,e'_4\}$) define the same orientation on $V$ (on
$V^{\perp}$, respectively). We define the function $\phi(t)$ by
    \bee\label{eq:case_iii_1_phi}
    \phi(t)=R(e'_1,\cos t e'_3+\sin t e'_4,e'_1,\cos t e'_3+\sin t
    e'_4).
    \eee
Then, we have $\phi'(0)=0$, and hence
    \be\label{eq:case_iii_1_1314}
    R'_{1314}=0.
    \ee
Similarly, considering the function $\psi(t)$ defined by
    \bee
    \psi(t)=R(\cos t e'_1+\sin t e'_2,e'_3,\cos t e'_1+\sin t
    e'_2,e'_3),
    \eee
we have $\psi'(0)=0$, and hence,
    \be\label{eq:case_iii_1_1323}
    R'_{1323}=0.
    \ee
Then, from \eqref{eq:case_iii_1_other}, \eqref{eq:case_iii_1_1314}
and \eqref{eq:case_iii_1_1323}, we see that the orthonormal basis
$\{e'_i\}$ is a generalized Singer-Thorpe basis. Similarly to
 Case (ii), we may choose a generalized
Singer-Thorpe basis for Cases (vi), (vii), (viii). Further, we may
also choose a generalized Singer-Thorpe basis for
 Cases III-2 and III-3.

\noindent{\bf Case IV-1.} Then, it suffices to consider Case (viii)
with $\lambda\ne0$. Then from \eqref{eq:i=1,j=2} and
\eqref{eq:i=1,j=3}, we have
    \be\label{eq:case_iv_1_other}
    \begin{gathered}
    R_{1223}=0,~R_{1213}=0,~R_{1323}=0.
    \end{gathered}
    \ee
Further, from Case (viii), we have
    \be\label{eq:case_vi+caseVIII}
    \bal
    &R_{1212}=R_{1313}=R_{2323}=-\lambda,\\
    &R_{1414}=R_{2424}=R_{3434}=\lambda.
    \eal
    \ee
Here, we note that each of the relations in
\eqref{eq:case_iv_1_other} and \eqref{eq:case_vi+caseVIII} is
preserved under the changes of the orthonormal basis satisfying the
conditions of Case IV-1. Let $V$ be a 3-dimensional subspace of
$T_pM$ spanned by $\{e_1,e_2,e_3\}$ satisfying that $V$ is
orthogonal complement of $\{e_4\}$.
    We define
    \be\label{eq:case_iv_1}
    R(e_1',e_2',e_2',e_4')={\underset{|x|=|y|=1}{\max_{x,~y\in
    V,~x\perp y}}}R(x,y,y,e_4),
    \ee where $e_3'\in V$ such that $e_3'\perp e_1'$, $e_3'\perp e_2'$, $|e_3'| =1,~ e_4'=e_4$.
First, we define the function $\phi(t)$ by
    \be\label{eq:case_vi_1_phi}
    \phi(t)=R(e_1', \cos t e_2'+\sin t e_3', \cos t e_2'+\sin t e_3',
    e_4').
    \ee
Then, by the hypothesis \eqref{eq:case_iv_1}, we have $\phi'(0)=0$,
and hence,
    \be\label{eq:case_iv_1_1234}
    R'_{1234}+R'_{1324}=0.
    \ee
Next, we consider the function $\psi(t)$ defined by
    \be\label{eq:case_iv_1_psi}
    \psi(t)=R(\cos t e_1' +\sin t e_3', e_2', e_2', e_4').
    \ee
Then we have $0=\psi'(0)=R'_{3224}$, and hence,
    \be\label{eq:case_iv_1_1314}
    R'_{1314}=0.
    \ee
Next, we consider the function $\zeta(t)$ defined by
    \be\label{eq:case_vi_1_zeta}
    \zeta(t)=R(\cos t e_1'+\sin t e_2',-\sin t e_1'+\cos t e_2',-\sin t
    e_1'+\cos t e_2', e_4').
    \ee
Then, by the hypothesis we have also $0=\zeta'(0)=-R'_{1214}$, and
hence
    \be\label{eq:case_iv_1_1214}
    R'_{1214}=0.
    \ee
Now, we set
    \be\label{eq:case_iv_1_basis}
    \begin{gathered}
    e''_2=\frac{1}{\sqrt{2}}e_2'+\frac{1}{\sqrt{2}}e_3',\\
    e''_3=-\frac{1}{\sqrt{2}}e_2'+\frac{1}{\sqrt{2}}e_3',\\
    e''_1=e_1',~e''_4=e_4'.
    \end{gathered}
    \ee
Then, we have
    \bee
    \bal
    R(e''_1,e''_2,e''_2,e''_4)=&\frac{1}{2}R(e_1',e_2'+e_3',e_2'+e_3',e_4')\\
    =&\frac{1}{2}\{R'_{1224}+R'_{1324}+R'_{1234}+R'_{1334}\}=0
    \eal
    \eee
by virtue of \eqref{eq:case_iv_1_1234}, and hence,
    \be\label{eq:case_iv_1_1224_o}
    R''_{1224}=0.
    \ee
Here, we set $R''_{ijkl}=R(e''_i,e''_j,e''_k,e''_l)$, $1\leq i,j,k,l
\leq 4$. Similarly, from \eqref{eq:case_iv_1_basis}, we have
    \be\label{eq:case_iv_1_1214_o}
    \bal
    R''_{1214}=\frac{1}{\sqrt{2}}R(e_1',e_2'+e_3',e_1',e_4')=\frac{1}{\sqrt{2}}(R'_{1214}+R'_{1314})=0,\\
    R''_{1314}=\frac{1}{\sqrt{2}}R(e_1',-e_2'+e_3',e_1',e_4')=\frac{1}{\sqrt{2}}(-R'_{1214}+R'_{1314})=0
    \eal
    \ee
by virtue of \eqref{eq:case_iv_1_1314} and
\eqref{eq:case_iv_1_1214}. Thus, from \eqref{eq:case_iv_1_1224_o}
and \eqref{eq:case_iv_1_1214_o}, we see that the orthonormal basis
$\{e''_i\}$ is a generalized Singer-Thorpe basis. Similarly, we may
also choose a generalized Singer-Thorpe basis for Cases
IV-2$\sim$IV-4.

Sufficiency: We assume that $M=(M,g)$ admits a generalized
Singer-Thorpe basis $\{e_i\}$. From the condition
\eqref{eq:general_S-T}, we see that \eqref{eq:i=1,j=2} and
\eqref{eq:i=1,j=3} hold on $M$. Further, by substituting
$\lambda_i=\sum_{k} R_{ikki}$ $(1\leq i\leq4)$ to the left hand
sides of \eqref{eq:i=j=1} and \eqref{eq:i=j=3}, and taking account
of \eqref{eq:from(5.1)}, we see also that each equation in
\eqref{eq:i=j=1} and \eqref{eq:i=j=3} holds. Therefore we see that
$M$ satisfies the curvature condition \eqref{eq:id+assumption}. Thus
$M$ is a weakly Einstein manifold by virtue of \eqref{eq:id}. This
completes the proof of Theorem B. \hfill$\square$\medskip

\section{An application}
In this section, we shall give a generalization of the Hitchin
inequality for a 4-dimensional compact oriented Einstein manifold.
Let $M=(M,g)$ be a compact oriented weakly Einstein manifold. Then,
from Theorem B, we may choose an generalized Singer-Thorpe basis
$\{e_i\}$ of $T_pM$ at any point $p\in M$ compatible with the
orientation of $M$. We set
    \be\label{eq:general_S-T_app}
    \begin{gathered}
    \alpha'_1=R_{1212},\qquad\alpha'_2=R_{1313},\qquad\alpha'_3=R_{1414},\\
    \alpha''_1=R_{3434},\qquad\alpha''_2=R_{2424},\qquad\alpha''_3=R_{2323},\\
    \beta_1=R_{1234},\qquad\beta_2=R_{1342},\qquad\beta_3=R_{1423}.
    \end{gathered}
    \ee
Then, from \eqref{eq:general_S-T_app}, by the first Bianchi
identity,
    \be
    \beta_1+\beta_2+\beta_3=0.
    \ee
Further, we set $ {\bf a}'=(\alpha'_1,\alpha'_2,\alpha'_3)$, $ {\bf
a}''=(\alpha''_1,\alpha''_2,\alpha''_3)$ and $ {\bf
b}=(\beta_1,\beta_2,\beta_3)$ and denote the canonical inner product
by $<,>$ on the 3-dimensional Euclidean space $\mathbb{R}^3$. We set
$| \bf{x}|=\sqrt{< \bf{x}, \bf{x}>}$ for any ${\bf
x}\in\mathbb{R}^3$. Then we may note that $|\bf{a}'|= |\bf{a}''|$ by
virtue of \eqref{eq:general-Singer}. Now, we denote the Euler number
and the first Pontrjagin number of $M$ by $\chi(M)$ and $p_{1}(M)$,
respectively. Then,
 from \eqref{eq:general_S-T_app}, applying the similar arguments in
\cite{Hi}, we have the following equalities:
    \be\label{eq:chi}
    \chi(M)=\frac{1}{4\pi^2}\int_M\{< {\bf a}', {\bf a}''>+| {\bf b}|^2\}dv_g
    \ee
and
    \be\label{eq:p_1}
    p_1(M)=\frac{1}{2\pi^2}\int_M< {\bf a}'+ {\bf a}'',
    {\bf b}>dv_g,
    \ee
where $dv_g$ is the volume element of $M$. Now, we set
    \be\label{eq:alpha}
    {\bf a}=\frac12({\bf a}'+{\bf a}'').
    \ee
Then, by \eqref{eq:alpha}, the equalities \eqref{eq:chi} and
\eqref{eq:p_1} are rewritten respectively by
    \be\label{eq:chi(M)}
    \chi(M)=\frac{1}{4\pi^2}\int_M\big\{2|{\bf a}|^2-|{\bf
    a}'|^2+|{\bf b}|^2\big\}dv_g,
    \ee
    \be\label{eq:p_1(M)}
    p_1(M)=\frac{1}{2\pi^2}\int_M2<{\bf a},{\bf b}>dv_g.
    \ee
Then, from \eqref{eq:chi(M)} and \eqref{eq:p_1(M)}, we have the
following:
    \be
    \bal
    &2\chi(M)\pm p_1(M)\\
    =&\frac{1}{2\pi^2}\int_M\big\{|{\bf a}|^2+|{\bf
    b}|^2\pm2<{\bf a},{\bf b}>+|{\bf a}|^2-|{\bf a}'|^2\big\}dv_g\\
    =&\frac{1}{2\pi^2}\int_M\big\{|{\bf a}\pm{\bf b}|^2+|{\bf a}|^2-|{\bf
    a}'|^2\big\}dv_g\\
    \geq&\frac{1}{2\pi^2}\int_M\big\{|{\bf a}|^2-|{\bf
    a}'|^2\big\}dv_g.
    \eal
    \ee
We set $f=|{\bf a}|^2-|{\bf a}'|^2$, Then, from the definition of
the vectors ${\bf a}'$, ${\bf a}''$ and ${\bf a}$, taking account of
the proof of Proposition \ref{th:4.1}, we have
    \be\label{eq:case(i)}
    f=0\qquad\qquad\qquad\qquad\qquad\qquad\quad\qquad\qquad\qquad\qquad\text{ for Case (i)},
    \ee
    \be\label{eq:case(ii)}
    f=-\frac{1}{4}(\lambda_1-\lambda_3)^2\quad(\lambda_1=\lambda_2,\lambda_3=\lambda_4)\qquad\qquad\qquad\text{
    for Case (ii)},
    \ee
    \be\label{eq:case(iii)}
    f=-\frac{1}{4}(\lambda_1-\lambda_2)^2\quad(\lambda_1=\lambda_3,\lambda_2=\lambda_4)\qquad\qquad\qquad\text{
    for Case (iii)},
    \ee
    \be\label{eq:case(iv)}
    f=-\frac{1}{4}(\lambda_1-\lambda_3)^2\quad(\lambda_1=\lambda_4,\lambda_2=\lambda_3)\qquad\qquad\qquad\text{
    for Case (iv)},
    \ee
    \be\label{eq:case(v)}
    f=-\frac{1}{4}\Big\{(\lambda_1-\lambda_3)^2+(\lambda_1-\lambda_4)^2\Big\}\quad(\lambda_1+\lambda_2=\lambda_3+\lambda_4)\quad\text{
    for Case (v)},
    \ee
    \be\label{eq:case(vi)}
    f=-\frac{1}{4}\Big\{(\lambda_1-\lambda_2)^2+(\lambda_1-\lambda_4)^2\Big\}\quad(\lambda_1+\lambda_3=\lambda_2+\lambda_4)\quad\text{
    for Case (vi)},
    \ee
    \be\label{eq:case(vii)}
    f=-\frac{1}{4}\Big\{(\lambda_1-\lambda_2)^2+(\lambda_1-\lambda_3)^2\Big\}\quad(\lambda_1+\lambda_4=\lambda_2+\lambda_3)\quad\text{
    for Case (vii)},
    \ee
    \be\label{eq:case(viii)}
    \bal
    f=&-\frac{1}{4}\Big\{(\lambda_1+\lambda_2)^2+(\lambda_1+\lambda_3)^2+(\lambda_1+\lambda_4)^2\Big\}\\
    &(\lambda_1+\lambda_2+\lambda_3+\lambda_4=0)\qquad\qquad\qquad\qquad\qquad\qquad\text{for Case (viii)}
    \eal
    \ee
at $p\in M$. Then from
\eqref{eq:case(ii)}$\sim$\eqref{eq:case(viii)}, we see that $f$
gives rise a continuous function on $M$ and further, $f=0$ holds at
$p$ if and only if $\lambda_1=\lambda_2=\lambda_3=\lambda_4$ holds
at $p$ (namely, $M$ is Einstein at $p$). Therefore, summing up the
above arguments we have finally the following Theorem.
    \vskip0.3cm

\noindent{\bf Theorem C} {\it Let $M=(M,g)$ be a compact weakly
Einstein manifold. Then, the following inequality holds on $M$:
    \be\label{eq:hit}
    2\chi(M)\pm p_1(M)\geq C,
    \ee
where $C=\frac{1}{2\pi^2}\int_M\{|{\bf a}|^2-|{\bf
a}'|^2\}dv_g\leq0$.} \vskip0.3cm
\begin{rem}\label{rm:4}
{\rm Since $p_1(M)=3\sigma(M)$ ($\sigma(M)$ is the Hirzebruch
signature of $M$), from Theorem C together with the proof, we see
that the inequality \eqref{eq:hit} reduces to the Hitchin inequality
\cite{Hi}
    \be\label{eq:hitchin}
    2\chi(M)\geq3|\sigma(M)|,
    \ee
    for the case where $M$ is Einstein. Thus,
    the inequality \eqref{eq:hit} in Theorem C is regarded as the generalization of the Hitchin
    inequality \eqref{eq:hitchin}.}
\end{rem}
The following example illustrates Theorem C and Remark \ref{rm:4}.
\begin{ex}\label{ex:6}
{\rm Let $M_1$ and $M_2$ be a unit 2-sphere and a compact oriented
surface of genus $m$ ($m\geq2$) with constant Gaussian curvature
$-1$, respectively, and further, $M$ be the Riemannian product of
$M_1$ and  $M_2$, $M=M_1\times M_2$. Then, we may easily check that
$M$ is a compact, oriented weakly Einstein manifold which is a
special case of Example \ref{ex:3}. Then, by taking account of the
K\"unneth formula, the Gauss-Bonnet formula and the formulas in
\cite{Hi}, we have
    \be\label{eq:ex6}
    \chi(M)=4(1-m),~p_1(M)=0 \text{ (thus,
    }\sigma(M)=0), \text{ and } C=8(1-m).
    \ee
Therefore, from \eqref{eq:ex6}, we see that the equality sign of the
inequality \eqref{eq:hit} in Theorem C holds for $M$, but $M$ does
not satisfy the Hitchin inequality \eqref{eq:hitchin}.}
\end{ex}

\section*{Acknowledgements}

 Research of Yunhee Euh was supported by the National
Research Foundation of Korea Grant funded by the Korean Government
[NRF-2009-352-C00007]. Research of JeongHyeong Park was supported by
Basic Science Research Program through the National Research
Foundation of Korea(NRF) funded by the Ministry of Education,
Science and Technology (2009-0087201).

\newcommand{\J}[4]{{\sl #1} {\bf #2} (#3) #4}
\newcommand{\CMP}{Comm.\ Math.\ Phys.}

\end{document}